\documentclass[reqno,11pt,a4paper]{amsart}
\usepackage[margin=2.7cm]{geometry}
\usepackage{cite}
\usepackage{newpxtext}
\usepackage{cabin}% sans serif
\usepackage[varqu,varl]{inconsolata}% sans serif typewriter
\usepackage[bigdelims,vvarbb]{newpxmath}% bb from STIX
\usepackage[cal=cm,bb=esstix,bbscaled=1.05]{mathalfa}% mathcal
\usepackage{mleftright}
\usepackage{graphicx}
\usepackage{colortbl}
\usepackage{booktabs}
\usepackage[up]{caption}
\usepackage[labelformat=simple]{subcaption}
\usepackage{xcolor} 
\usepackage{hyperref}
\hypersetup{
    colorlinks,
    linkcolor={red!50!black},
    citecolor={blue!50!black},
    urlcolor={blue!80!black}
}
\AtBeginDocument{%
\def\MR#1{}
}
\newcommand{\orcid}[1]{\,\resizebox{8px}{!}{\href{https://orcid.org/#1}{\includegraphics{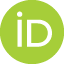}}}}
\graphicspath{{images_300/}} % 300 dpi
\newcommand{\ZZ}{\mathbb{Z}}

\newcommand{\QQ}{\mathbb{Q}}
\newcommand{\RR}{\mathbb{R}}

\newcommand{\PP}{\mathbb{P}}
\newcommand{\NQ}{N_\QQ}
\newcommand{\MQ}{M_\QQ}

\DeclareMathOperator{\GL}{GL}
\DeclareMathOperator{\conv}{conv}
\DeclareMathOperator{\Ehr}{Ehr}
\DeclareMathOperator{\Vol}{Vol}
\DeclareMathOperator{\vol}{vol}
\DeclareMathOperator{\Hilb}{Hilb}
\DeclareMathOperator{\Hom}{Hom}

\theoremstyle{plain}
\newtheorem{Lemma}{Lemma}[section]

\newtheorem{question}{Question}
\newtheorem{example}[Lemma]{Example}
\theoremstyle{remark}
\newtheorem{algorithm}[Lemma]{Algorithm}
\numberwithin{equation}{section}
\begin{document}
%-------------------------------------------------------------------------------
% Address details for each author
\author[T. Coates]{Tom Coates\orcid{0000-0003-0779-9735}}
\address{Department of Mathematics\\Imperial College London\\180 Queen's Gate\\London\\SW7 2AZ\\UK}
\email{t.coates@imperial.ac.uk}% Coates
\author[J.\,Hofscheier]{Johannes Hofscheier\orcid{0000-0001-8642-3984}}
\address{School of Mathematical Sciences\\University of Nottingham\\Nottingham\\NG7 2RD\\UK}
\email{johannes.hofscheier@nottingham.ac.uk}% Hofscheier
\author[A.\,M.\,Kasprzyk]{Alexander M.~Kasprzyk\orcid{0000-0003-2340-5257}}
\address{School of Mathematical Sciences\\University of Nottingham\\Nottingham\\NG7 2RD\\UK}
\email{a.m.kasprzyk@nottingham.ac.uk}% Kasprzyk
%-------------------------------------------------------------------------------
%-------------------------------------------------------------------------------
\title{Machine Learning the Dimension of a Polytope}
%-------------------------------------------------------------------------------
\begin{abstract}
We use machine learning to predict the dimension of a lattice polytope directly from its Ehrhart series. This is highly effective, achieving almost~$100\%$ accuracy. We also use machine learning to recover the volume of a lattice polytope from its Ehrhart series, and to recover the dimension, volume, and quasi-period of a rational polytope from its Ehrhart series. In each case we achieve very high accuracy, and we propose mathematical explanations for why this should be so.
\end{abstract}
%-------------------------------------------------------------------------------
\maketitle
%-------------------------------------------------------------------------------
\section{Introduction}\label{sec:intro}
%-------------------------------------------------------------------------------
Let~$P\subset\ZZ^d\otimes_\ZZ\QQ$ be a convex lattice polytope of dimension~$d$, that is, let~$P$ be the convex hull of finitely many points in~$\ZZ^d$ whose~$\QQ$-affine combinations generate~$\QQ^d$. A fundamental invariant of~$P$ is the number of lattice points that it contains,~$|P\cap\ZZ^d|$. More generally, let~$L_P(m)\coloneq|mP\cap\ZZ^d|$ count the lattice points in the~$m$th dilation~$mP$ of~$P$, where~$m\in\ZZ_{\geq 0}$. Then~$L_P$ is given by a polynomial of degree~$d$ called the~\emph{Ehrhart polynomial}~\cite{Ehr62}. The corresponding generating series, called the~\emph{Ehrhart series} and denoted by~$\Ehr_P$, can be expressed as a rational function with numerator a polynomial of degree at most~$d$~\cite{Stan80}:
\begin{align*}
\Ehr_P(t)\coloneq\sum_{m\geq 0}L_P(m)t^m=\frac{\delta_0+\delta_1t+\cdots+\delta_dt^d}{(1-t)^{d+1}}, && \delta_i \in \ZZ.
\end{align*}
The coefficients~$(\delta_0,\delta_1,\ldots,\delta_d)$ of this numerator, called the~\emph{$\delta$-vector} or~\emph{$h^*$-vector} of~$P$, have combinatorial meaning~\cite{Ehr67}:
\begin{enumerate}
\item $\delta_0=1$;
\item $\delta_1=|P\cap\ZZ^d|-d-1$;
\item $\delta_d=|P^\circ\cap\ZZ^d|$, where~$P^\circ=P\setminus\partial P$ is the strict interior of~$P$;
\item $\delta_0+\cdots+\delta_d=\Vol(P)$, where~$\Vol(P)=d!\vol(P)$ is the lattice-normalised volume of~$P$.
\end{enumerate}
The polynomial~$L_P$ can be expressed in terms of the~$\delta$-vector:
\[
L_P(m)=\sum_{i=0}^d\delta_i\binom{d+m-i}{d}.
\]
From this we can see that the leading coefficient of~$L_P$ is~$\vol(P)$, the Euclidean volume of~$P$.

Given~$d+1$ terms of the Ehrhart series, one can recover the~$\delta$-vector and hence the invariants~$|P\cap\ZZ^d|$,~$|P^\circ\cap\ZZ^d|$, and~$\Vol(P)$. This, however, assumes knowledge of the dimension~$d$.

\begin{question}\label{qu:ML_dimension}
Given a lattice polytope~$P$, can machine learning recover the dimension~$d$ of~$P$ from sufficiently many terms of the Ehrhart series~$\Ehr_P$?
\end{question}

There has been recent success using machine learning~(ML) to predict invariants such as~$\Vol(P)$ directly from the vertices of~$P$~\cite{BHHH+21}, and to predict numerical invariants from a geometric analogue of the Ehrhart series called the Hilbert series~\cite{BHHH+22}. As we will see in~\S\S\ref{sec:lattice_data}--\ref{sec:volume}, ML is also extremely effective at answering Question~\ref{qu:ML_dimension}. In~\S\ref{sec:log_asymptotics} we propose a possible explanation for this.
%-------------------------------------------------------------------------------
\subsection{The quasi-period of a rational convex polytope}\label{sec:quasiperiod}
%-------------------------------------------------------------------------------
Now let~$P$ be a convex polytope with~\emph{rational} vertices, and let~$k\in\ZZ_{>0}$ be the smallest positive dilation of~$P$ such that~$kP$ is a lattice polytope. One can define the Ehrhart series of~$P$ exactly as before:
\[
	L_P(m)\coloneq|mP\cap\ZZ^d|.
\]
In general~$L_P$ will no longer be a polynomial, but it is a~\emph{quasi-polynomial} of degree~$d$ and period~$k$~\cite{Ehr62,McMul78}. That is, there exist polynomials~$f_0,\dots, f_{k-1}$, each of degree~$d$, such that:
\begin{align*}
	L_P(q k + r) = f_r(q),
	&& \text{whenever~$0 \leq r < k$.}
\end{align*}
The leading coefficient of each~$f_r$ is~$\vol(kP)$, the Euclidean volume of the~$k$th dilation of~$P$. 

It is sometimes possible to express~$L_P$ using a smaller number of polynomials. Let~$\rho$ be the minimum number of polynomials needed to express~$L_P$ as a quasi-polynomial;~$\rho$ is called the~\emph{quasi-period} of~$P$, and is a divisor of~$k$. Each of these~$\rho$ polynomials is of degree~$d$ with leading coefficient~$\vol(\rho P)$. When~$\rho<k$ we say that~$P$ exhibits~\emph{quasi-period collapse}~\cite{MW05,BSW08,FK08,KW18,BER21}. The Ehrhart series of~$P$ can be expressed as a rational function:
\begin{align*}
\Ehr_P(t)=\frac{\delta_0+\delta_1t+\cdots+\delta_{\rho(d+1)-1}t^{\rho(d+1)-1}}{(1-t^\rho)^{d+1}}, && \delta_i \in \ZZ.
\end{align*}
As in the case of lattice polytopes, the~$\delta$-vector carries combinatorial information about~$P$~\cite{McA09,Her10,BER21,BBVM22,HHK}. Given~$\rho(d+1)$ terms of the Ehrhart series, one can recover the~$\delta$-vector. This, however, requires knowing both the dimension~$d$ and the quasi-period~$\rho$ of~$P$.

\begin{question}\label{qu:ML_quasiperiod}
Given a rational polytope~$P$, can machine learning recover the dimension~$d$ and quasi-period~$\rho$ of~$P$ from  sufficiently many terms of the Ehrhart series~$\Ehr_P$?
\end{question}

\noindent ML again performs very well here -- see~\S\ref{sec:ML_quasiperiod}. We propose a mathematical explanation for this, in terms of forward differences and affine hyperplanes, in~\S\ref{sec:forward_differences}. 
%-------------------------------------------------------------------------------
\subsection{Some motivating examples}
%-------------------------------------------------------------------------------
\begin{example}\label{eg:P2}
	The~$2$-dimensional lattice polytope~$P\coloneq\conv\{(-1, -1),(-1, 2),(2, -1)\}$ has volume~$\Vol(P)=9$,~$|P\cap\ZZ^2|=10$, and~$|P^\circ\cap\ZZ^2|=1$. The Ehrhart polynomial of~$P$ is:
	\[
	L_P(m)=\frac{9}{2}m^2 + \frac{9}{2}m + 1,
	\]
	and the Ehrhart series of~$P$ is generated by:
	\[
	\Ehr_P(t)=\frac{1+7t+t^2}{(1-t)^3}=1+10t+28t^2+55t^3+91t^4+\cdots.
	\]
\end{example}

\begin{example}\label{eg:P114}
	The smallest dilation of the triangle~$P\coloneq\conv\{(5,-1),(-1,-1),(-1,1/2)\}$ giving a lattice triangle is~$2P$. Hence~$L_P$ can be written as a quasi-polynomial of degree~$2$ and period~$2$:
	\[
	f_0(q)=18q^2 + 9q + 1,\qquad f_1(q)=18q^2 + 27q + 10.
	\]
	This is not, however, the minimum possible period. As in Example~\ref{eg:P2},
	\[
	L_P(m)=\frac{9}{2}m^2 + \frac{9}{2}m + 1
	\]
	and thus~$P$ has quasi-period~$\rho=1$. 
\end{example}

\noindent This striking example of quasi-period collapse is developed further in Example~\ref{eg:mutation_P2}.
%-------------------------------------------------------------------------------
\subsection{Code and data availability}
%-------------------------------------------------------------------------------
The datasets used in this work were generated with~V2.25-4 of Magma~\cite{BCP97}. We performed our~ML analysis with scikit-learn~\cite{scikit-learn}, a standard machine learning library for Python, using scikit-learn~v0.24.1 and Python~v3.8.8. All data, along with the code used to generate it and to perform the subsequent analysis, is available from Zenodo~\cite{dimData,quasiperiodData} under a permissive open source license (MIT for the code, CC0 for the data).
%-------------------------------------------------------------------------------
\section{Question~\ref{qu:ML_dimension}: Dimension}\label{sec:ML_dimension} 
%-------------------------------------------------------------------------------
In this section we investigate whether machine learning can predict the dimension~$d$ of a lattice polytope~$P$ from sufficiently many terms of the Ehrhart series~$\Ehr_P$. We calculate terms~$L_P(m)$ of the Ehrhart series, for~$0\leq m\leq 1100$, and encode these in a~\emph{logarithmic Ehrhart vector}:
\[
(\log y_0,\log y_1,\ldots,\log y_{1100}),\qquad\text{ where }y_m\coloneq L_P(m).
\]
We find that standard~ML techniques are extremely effective, recovering the dimension of~$P$ with almost~$100\%$ accuracy from the logarithmic Ehrhart vector. We then ask whether~ML can recover~$\Vol(P)$ from this Ehrhart data; again, this is achieved with near~$100\%$ success.
%-------------------------------------------------------------------------------
\subsection{Data generation}\label{sec:lattice_data}
%-------------------------------------------------------------------------------
A dataset~\cite{dimData} containing~$2918$ distinct entries, with~$2\leq d\leq 8$, was generated using Algorithm~\ref{alg:dimension} below. The distribution of this data is summarised in Table~\ref{tab:dimension}.

\begin{algorithm}\label{alg:dimension}\phantom{x}
\begin{description}
\item[Input] A positive integer~$d$.
\item[Output] A vector
\[
\big(\log y_0, \log y_1, \ldots, \log y_{1100}, d,\Vol(P)\big)
\]
for a~$d$-dimensional lattice polytope~$P$, where~$y_m\coloneq L_P(m)$.
\end{description}
\begin{enumerate}
\item\label{step:gen1}
Choose~$d+k$ lattice points~$\{v_1,\ldots,v_{d+k}\}$ uniformly at random in a box~$[-5,5]^d$, where~$k$ is chosen uniformly at random in~$\{1,\ldots,5\}$.
\item
Set~$P\coloneq\conv\{v_1,\ldots,v_{d+k}\}$. If~$\dim(P)\neq d$ return to step~\eqref{step:gen1}.
\item
Calculate the coefficients~$y_m\coloneq L_P(m)$ of the Ehrhart series of~$P$, for~$0\leq m\leq 1100$.
\item
Return the vector~$\big(\log y_0, \log y_1, \ldots, \log y_{1100}, d,\Vol(P)\big)$.
\end{enumerate}
\end{algorithm}

\noindent We deduplicated on the vector~$\big(\log y_0, \log y_1, \ldots, \log y_{1100}, d,\Vol(P)\big)$ to get a dataset with distinct entries. In particular, two polytopes that are equivalent under the group of affine-linear transformations~$\GL_d(\ZZ) \ltimes \ZZ^d$ give rise to the same point in the dataset.
%-------------------------------------------------------------------------------
\subsection{Machine learning the dimension}\label{sec:dimension}
\begin{figure}[tb]
\includegraphics[width=0.7\textwidth]{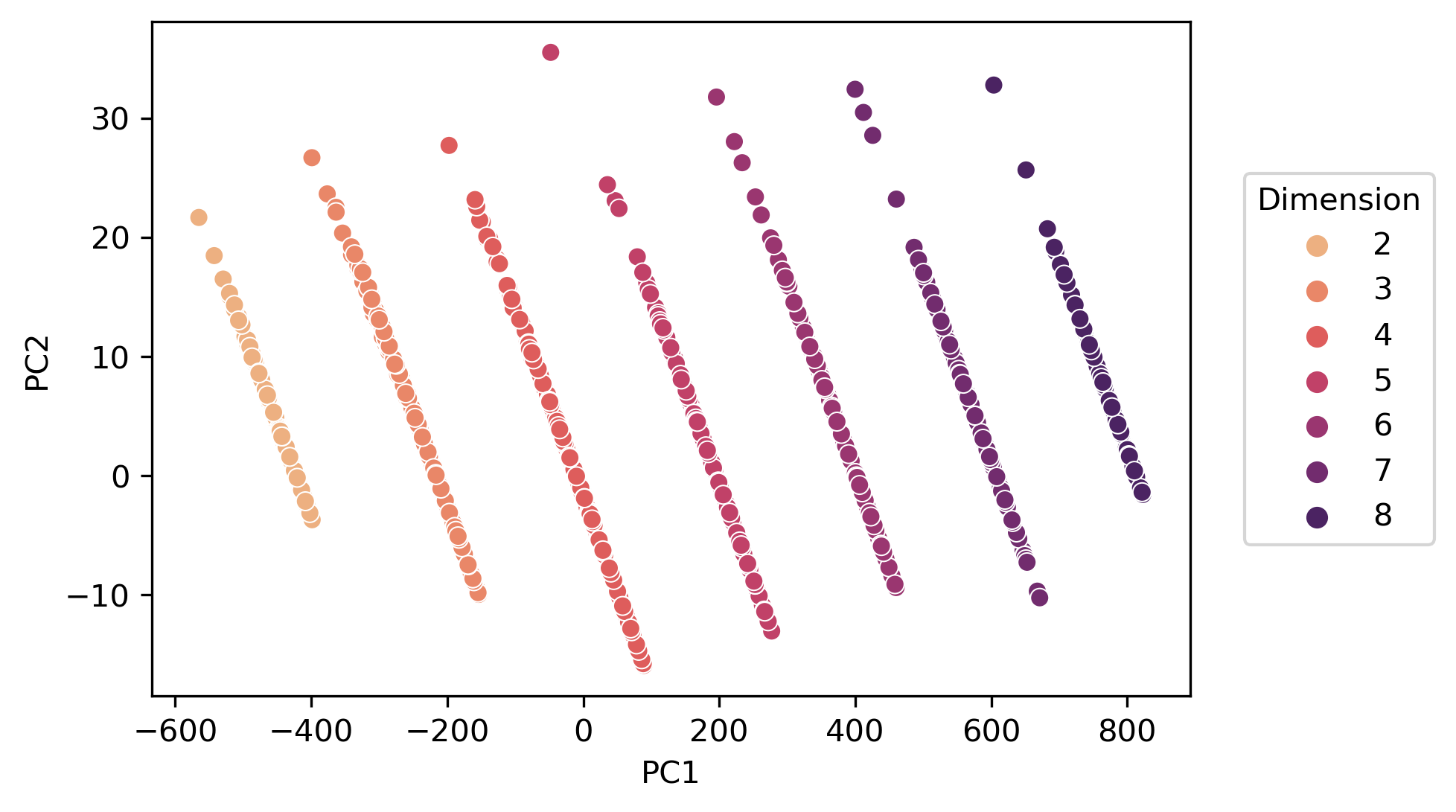}
\caption{Projection onto the first two principal components of the logarithmic Ehrhart vector coloured by the dimension of~$P$.\label{fig:pca_dim}}
\end{figure}

%-------------------------------------------------------------------------------
We reduced the dimensionality of the dataset by projecting onto the first two principal components of the logarithmic Ehrhart vector. As one would expect from Figure~\ref{fig:pca_dim}, a linear support vector machine (SVM) classifier trained on these features predicted the dimension of~$P$ with~$100\%$ accuracy. Here we used a scikit-learn pipeline consisting of a~\texttt{StandardScaler} followed by an~\texttt{SVC} classifier with linear kernel and regularization hyperparameter~$C=0.1$, using~$50\%$ of the data for training the classifier and tuning the hyperparameter, and holding out the remaining~$50\%$ of the data for model validation.

In~\S\ref{sec:log_asymptotics} below we give a mathematical explanation for the structure observed in Figure~\ref{fig:pca_dim}, and hence for why~ML is so effective at predicting the dimension of~$P$. Note that the discussion in~\S\ref{sec:forward_differences} suggests that one should also be able to extract the dimension using~ML on the~\emph{Ehrhart vector}
\[
	(y_0,y_1,\ldots,y_{1100})
\]
rather than the logarithmic Ehrhart vector
\[
	(\log y_0,\log y_1,\ldots,\log y_{1100}).
\]
This is indeed the case, although here it is important not to reduce the dimensionality of the data too much\footnote{This is consistent with the discussion in~\S\ref{sec:forward_differences}, which suggests that we should try to detect whether the Ehrhart vector lies in a union of linear subspaces that have fairly high codimension.}. A linear~SVM classifier trained on the full Ehrhart vector predicts the dimension of~$P$ with~$98.7\%$ accuracy, with a pipeline exactly as above except that~$C = 50000$. A linear~SVM classifier trained on the first~$30$ principal components of the Ehrhart vector (the same pipeline, but with~$C = 20$) gives~$93.7\%$ accuracy, but projecting to the first two components reduces accuracy to~$43.3\%$.

\begin{table}[tb]
	\small
	\centering
	\begin{tabular}{rccccccc}
	\toprule
	Dimension&2&3&4&5&6&7&8\\
	\midrule
	Total&431&787&812&399&181&195&113\\
	\bottomrule
	\end{tabular}
	\caption{The distribution of the dimensions appearing in the dataset for Question~\ref{qu:ML_dimension}.}
	\label{tab:dimension}
\end{table}

%-------------------------------------------------------------------------------
\subsection{Machine learning the volume}\label{sec:volume}
\begin{figure}[tb]
	\includegraphics[width=.45\textwidth]{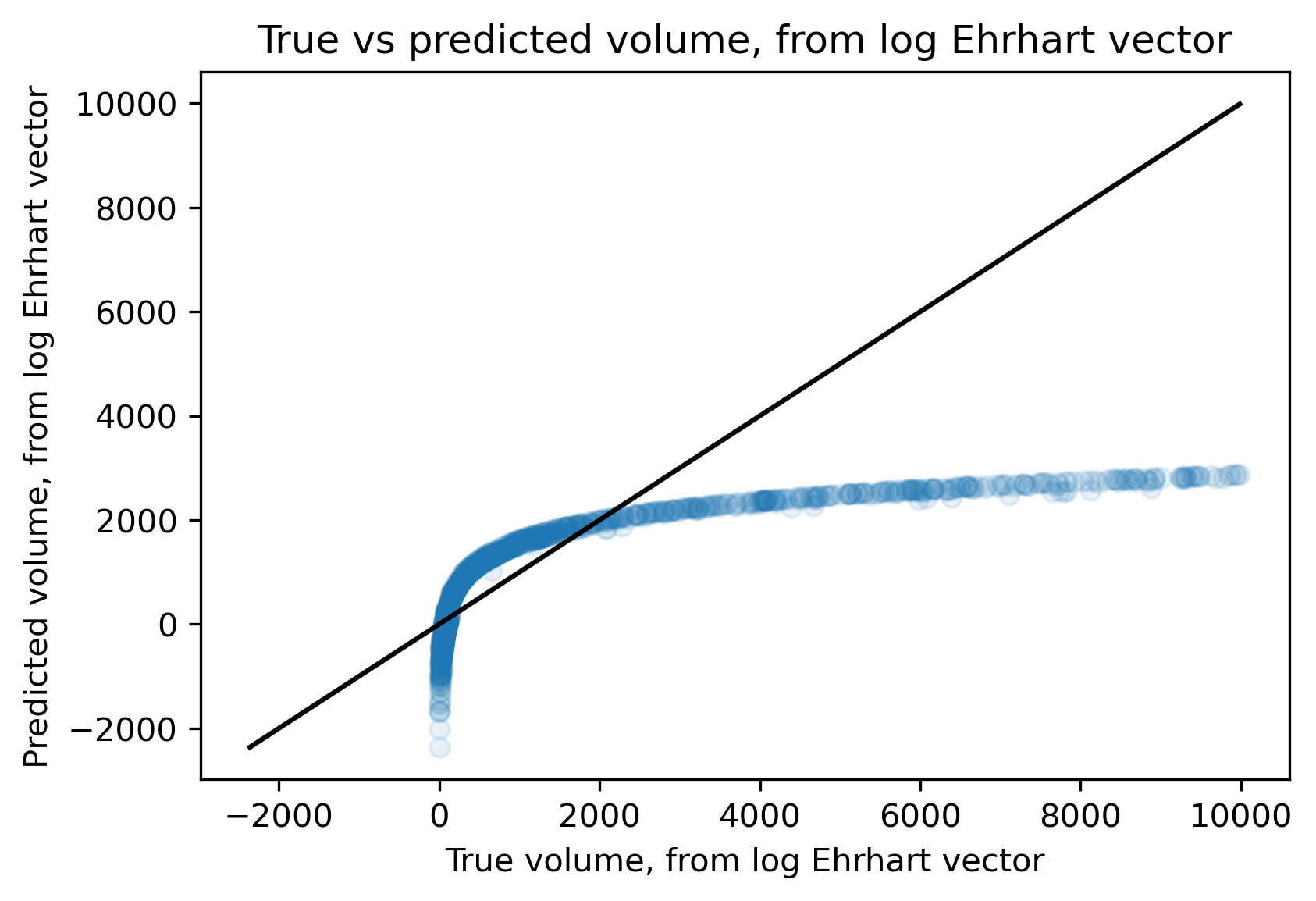}
	\includegraphics[width=.45\textwidth]{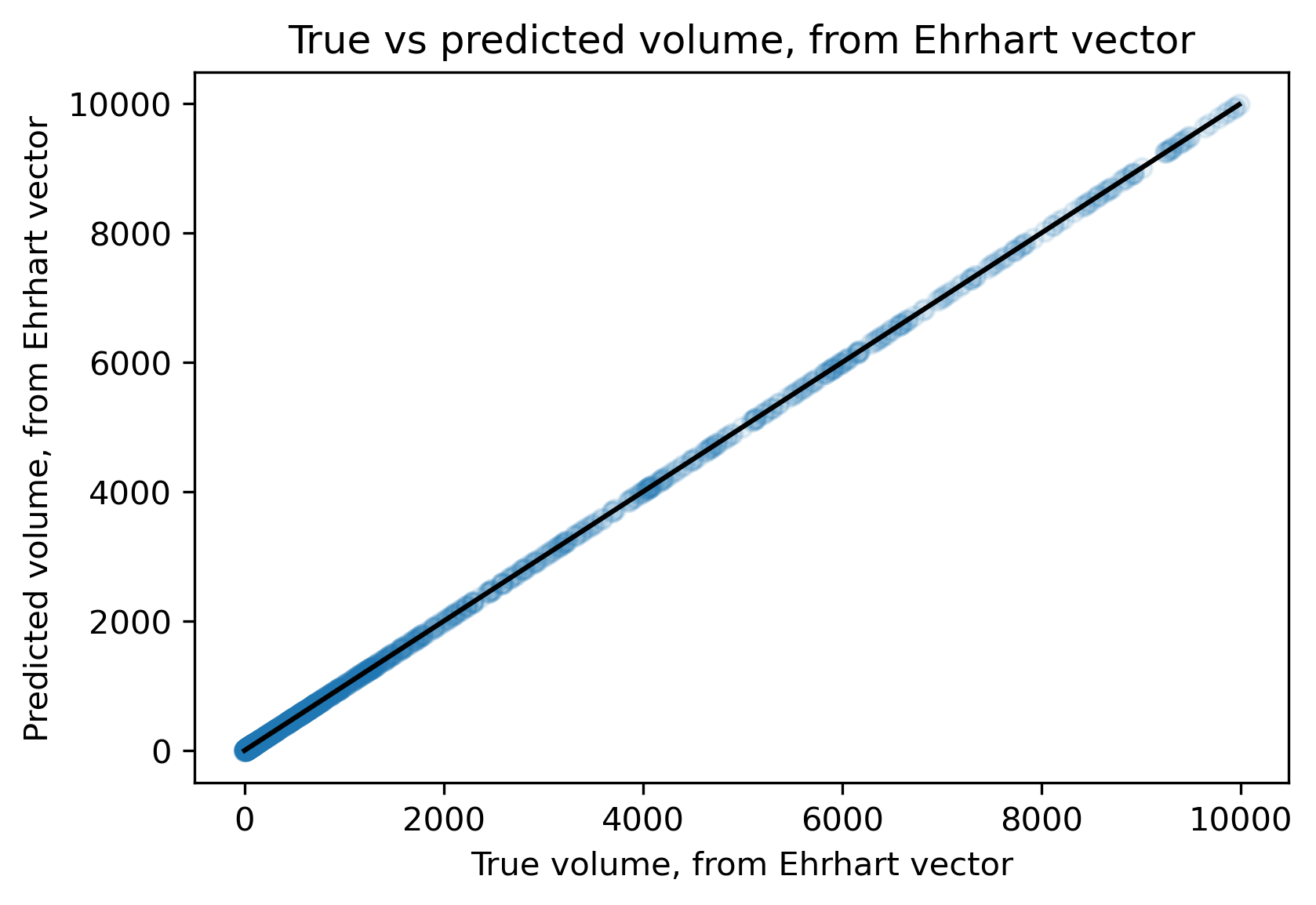}
	\caption{Plotting the true versus the predicted normalised volume.\label{fig:true_vs_predicted_volume}}
\end{figure}

%-------------------------------------------------------------------------------
To learn the normalised volume of a lattice polytope from its logarithmic Ehrhart vector, we used a scikit-learn pipeline consisting of a~\texttt{StandardScaler} followed by an~\texttt{SVR} regressor with linear kernel and regularization hyperparameter~$C = 340$. We restricted attention to the roughly~$75\%$ of the data with volume less than~$10\,000$, thereby removing outliers. We used~$50\%$ of that data for training and hyperparameter tuning, selecting the training set using a shuffle stratified by volume; this corrects for the fact that the dataset contains a high proportion of polytopes with small volume. The regression had a coefficient of determination ($R^2$) of~$0.432$, and gave a strong hint (see Figure~\ref{fig:true_vs_predicted_volume}) that we should repeat the analysis replacing the logarithmic Ehrhart vector with the Ehrhart vector.

Using the same pipeline (but with~$C=1000$) and the Ehrhart vector gives a regression with~$R^2=1.000$; see Figure~\ref{fig:true_vs_predicted_volume}. This regressor performs well over the full dataset, with volumes ranging up to approximately $4.5$~million: over the full dataset we still find~$R^2=1.000$. The fact that Support Vector Machine methods are so successful in recovering the volume of~$P$ from the Ehrhart vector is consistent with the discussion in~\S\ref{sec:forward_differences}.

%-------------------------------------------------------------------------------
\subsection{\texorpdfstring{\boldmath Crude asymptotics for~$\log y_k$}{Crude asymptotics for~log y\_k}}\label{sec:log_asymptotics}
\begin{figure}[tb]
	\includegraphics[width=0.7\textwidth]{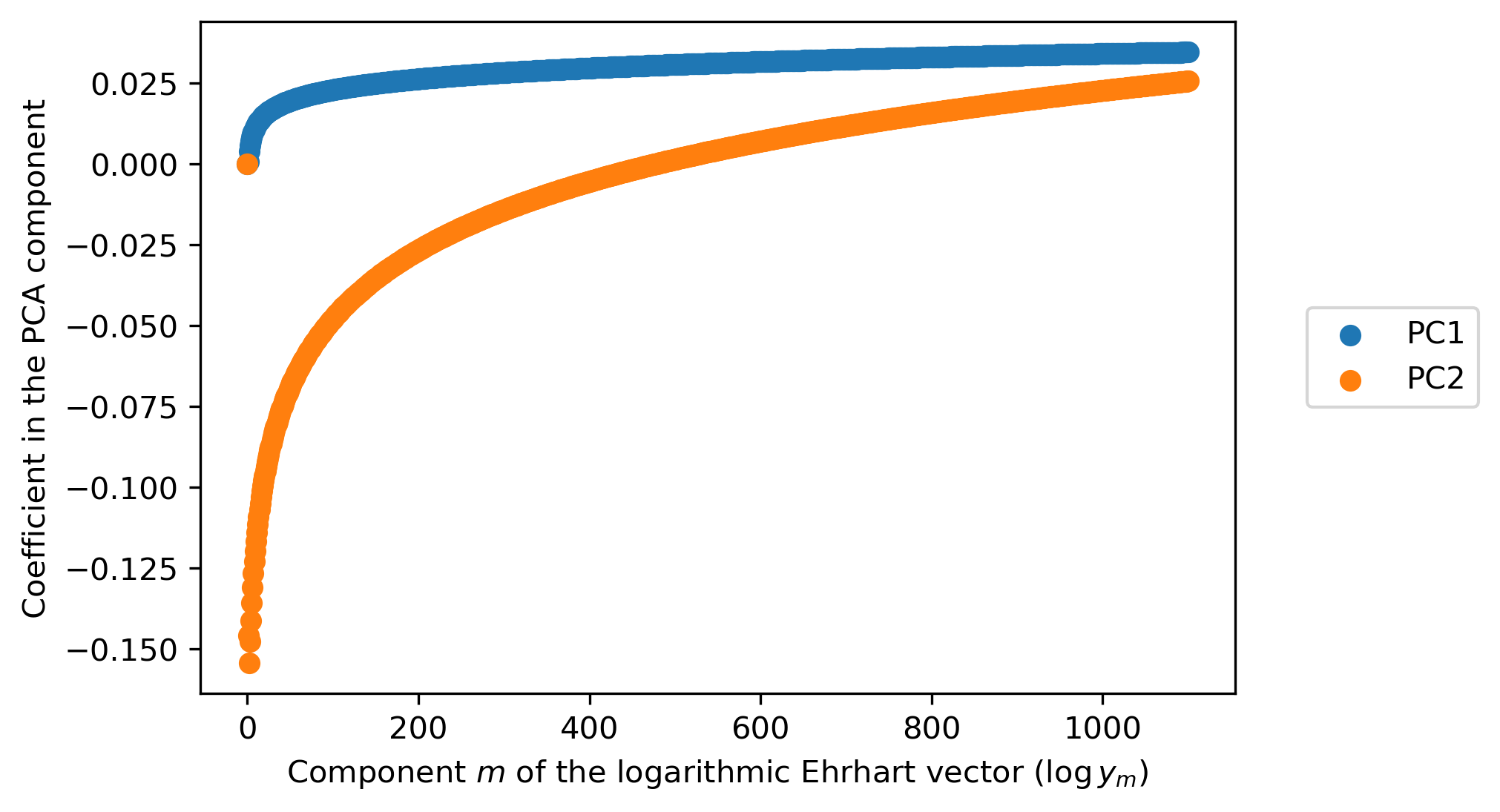}
	\caption{Contribution of~$\log y_m$ to the first two principal components of the logarithmic Ehrhart vector, as~$m$ varies.
	\label{fig:pca_dependence}}
\end{figure}

%-------------------------------------------------------------------------------
Since
\[
	L_P(m) = \vol(P) m^d + \text{lower order terms in~$m$}
\]
we have that 
\[
	\log y_m \sim d \log m + \log \vol(P).
\] 
For~$m \gg 0$, therefore, we see that the different components~$\log y_m$ of the logarithmic Ehrhart vector depend approximately affine-linearly on each other as~$P$ varies. It seems intuitively plausible that the first two~PCA components of the logarithmic Ehrhart vector should depend non-trivially on~$\log y_m$ for~$m \gg 0$, and in fact this is the case -- see Figure~\ref{fig:pca_dependence}. Thus the first two~PCA components of the logarithmic Ehrhart vector should vary approximately affine-linearly as~$P$ varies, with constant slope and with a translation that depends only on the dimension of~$P$.

%-------------------------------------------------------------------------------
\section{Question~\ref{qu:ML_quasiperiod}: Quasi-period}\label{sec:ML_quasiperiod}
%-------------------------------------------------------------------------------
Here we investigate whether~ML can predict the quasi-period~$\rho$ of a~$d$-dimensional rational polytope~$P$ from sufficiently many terms of the Ehrhart vector of~$P$.
Once again, standard~ML techniques based on Support Vector Machines are highly effective, achieving classification accuracies of up to~$95.3\%$. We propose a potential mathematical explanation for this in \S\ref{sec:forward_differences}.
%-------------------------------------------------------------------------------
\subsection{Data generation}
%-------------------------------------------------------------------------------
In this section we consider a dataset~\cite{quasiperiodData} containing~$84\,000$ distinct entries. This was obtained by using Algorithm~\ref{alg:quasiperiod} below to generate a larger dataset, followed by random downsampling to a subset with 2000 datapoints for each pair~$(d, \rho)$ with~$d \in \{2,3,4\}$ and~$\rho \in \{2,3,\ldots,15\}$.

\begin{algorithm}\label{alg:quasiperiod}\phantom{x}
\begin{description}
\item[Input] A positive integer~$d$.
\item[Output] A vector
\[
\big(\log y_0, \log y_1, \ldots, \log y_{1100}, d, \rho\big)
\]
for a~$d$-dimensional rational polytope~$P$ with quasi-period~$\rho$, where~$y_m\coloneq L_P(m)$.
\end{description}
\begin{enumerate}
\item Choose~$r\in \{2,3,\ldots,15\}$ uniformly at random.
\item\label{step:gen2} Choose~$d+k$ lattice points~$\{v_1,\ldots,v_{d+k}\}$ uniformly at random in a box~$[-5r,5r]^d$, where~$k$ is chosen uniformly at random in~$\{1,\ldots,5\}$.
\item
Set~$P\coloneq\conv\{v_1,\ldots,v_{d+k}\}$.  If~$\dim(P)\neq d$ return to step~\eqref{step:gen2}.
\item
Choose a lattice point~$v\in P\cap\ZZ^d$ uniformly at random and replace~$P$ with the translation~$P-v$. (We perform this step to ensure that the resulting rational polytope always contains a lattice point; this avoids complications when taking~$\log$ in step~\eqref{step:log}.)
\item
Replace~$P$ with the dilation~$P/r$.
\item
Calculate the coefficients~$y_m\coloneq L_P(m)$ of the Ehrhart series of~$P$, for~$0\leq m\leq 1100$.
\item Calculate the quasi-period~$\rho$.
\item\label{step:log}
Return the vector~$\big(\log y_0, \log y_1, \ldots, \log y_{1100}, d, \rho\big)$.
\end{enumerate}
\end{algorithm}
\noindent As before, we deduplicated the dataset on the vector~$\big(\log y_0, \log y_1, \ldots, \log y_{1100}, d, \rho\big)$.
%-------------------------------------------------------------------------------
\subsection{Recovering the dimension and volume}
\begin{figure}[b]
	\includegraphics[width=0.7\textwidth]{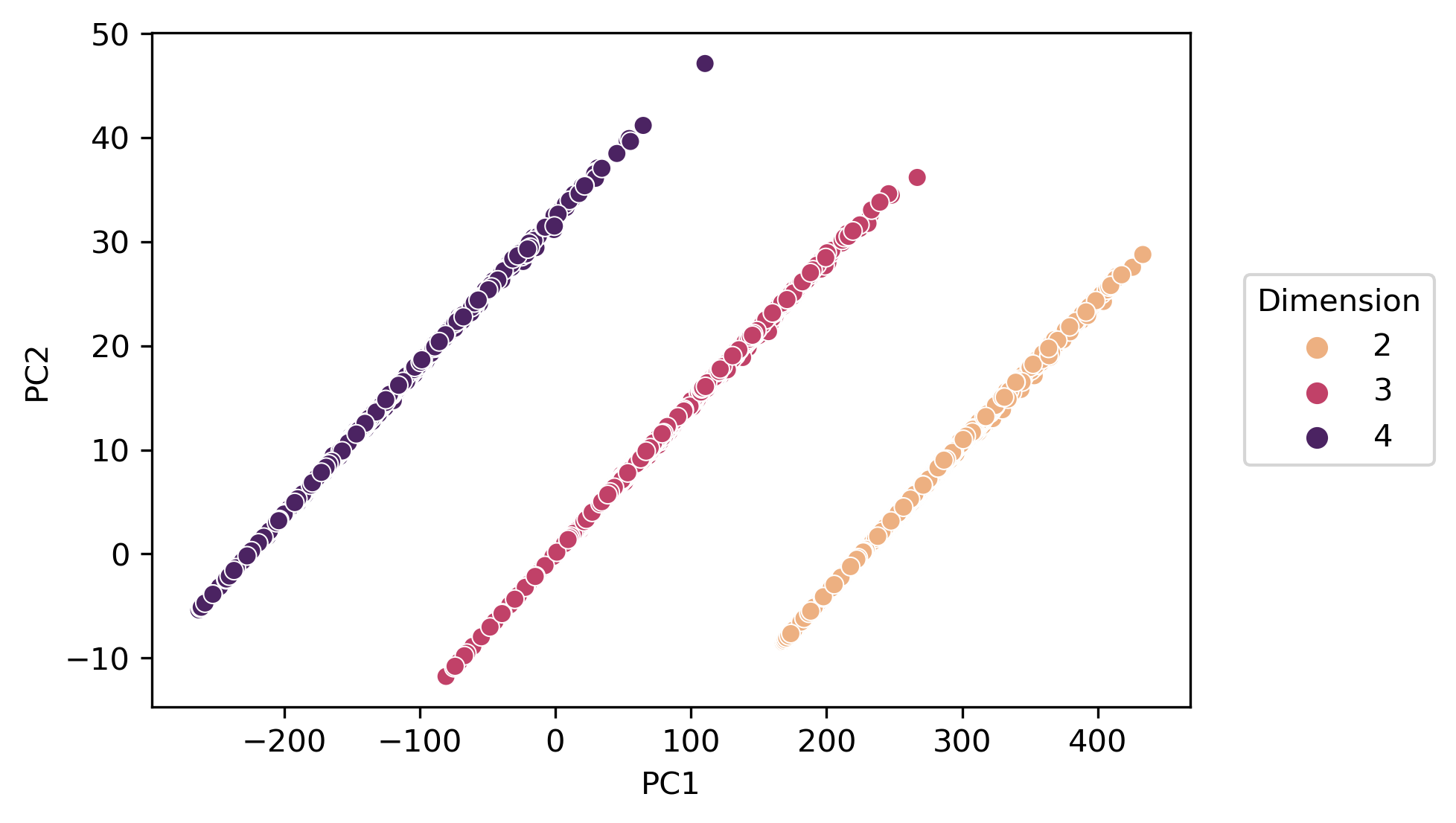}
	\caption{Projection onto the first two principal components of the logarithmic Ehrhart vector data coloured by the dimension of~$P$.\label{fig:pca_dim_rational}}
\end{figure}

%-------------------------------------------------------------------------------
Figure~\ref{fig:pca_dim_rational} shows the first two principal components of the logarithmic Ehrhart vector. As in~\S\ref{sec:dimension}, this falls into widely-separated linear clusters according to the value of~$\dim(P)$, and so the dimension of the rational polytope~$P$ can be recovered with high accuracy from its logarithmic Ehrhart vector. Furthermore, as in~\S\ref{sec:volume}, applying a linear~\texttt{SVR} regressor (with~$C=1000$) to the Ehrhart vector predicts the volume of a rational polytope~$P$ with high accuracy ($R^2=1.000$).

%-------------------------------------------------------------------------------
\subsection{Machine learning the quasi-period}
\begin{figure}[tb]
% Note: using "width" here doesn't work, because the rhs plots are wider. Using "heigh" gives
% the best result, even if it will need redoing for different layouts.
	\includegraphics[height=5.5cm]{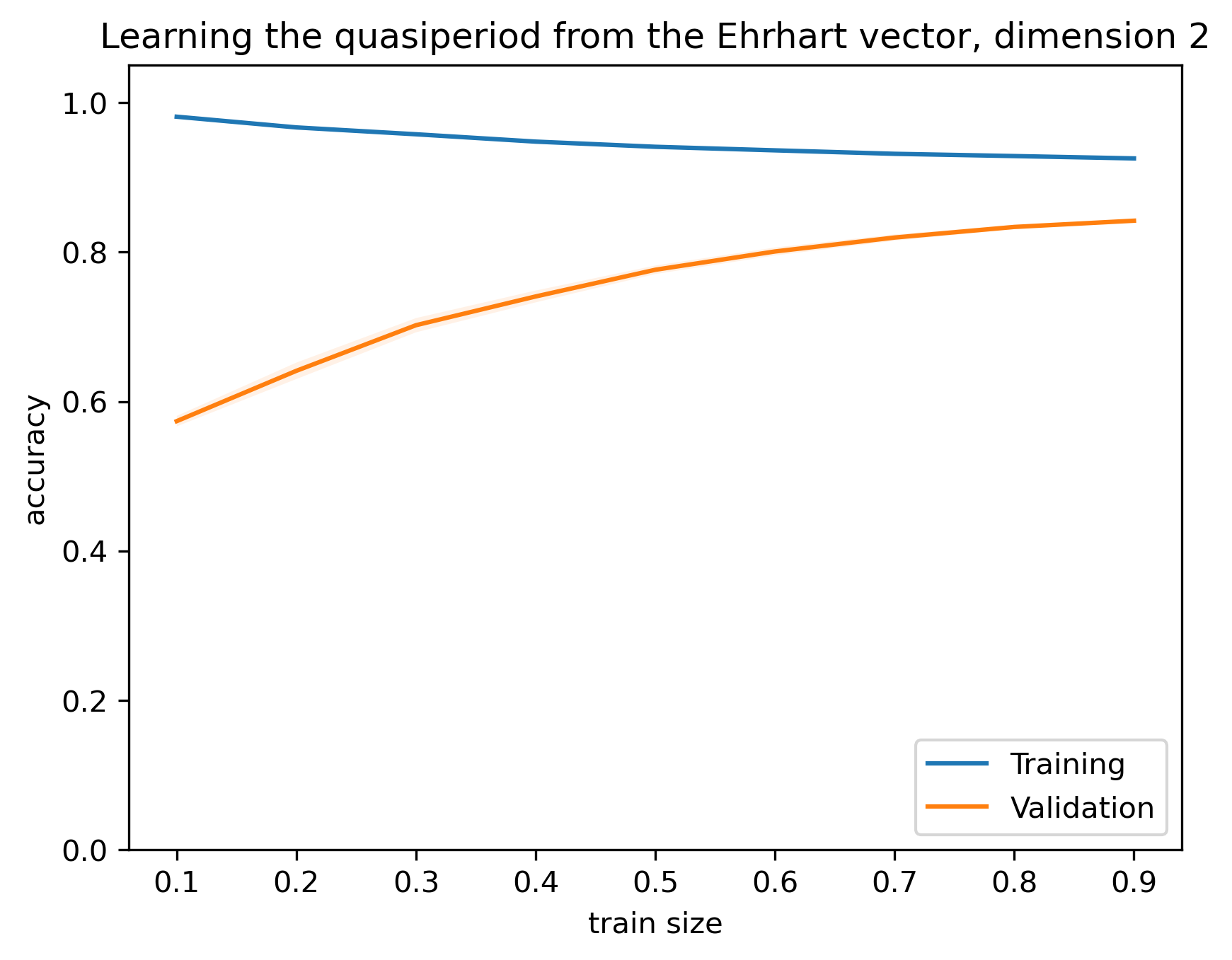}
	\includegraphics[height=5.5cm]{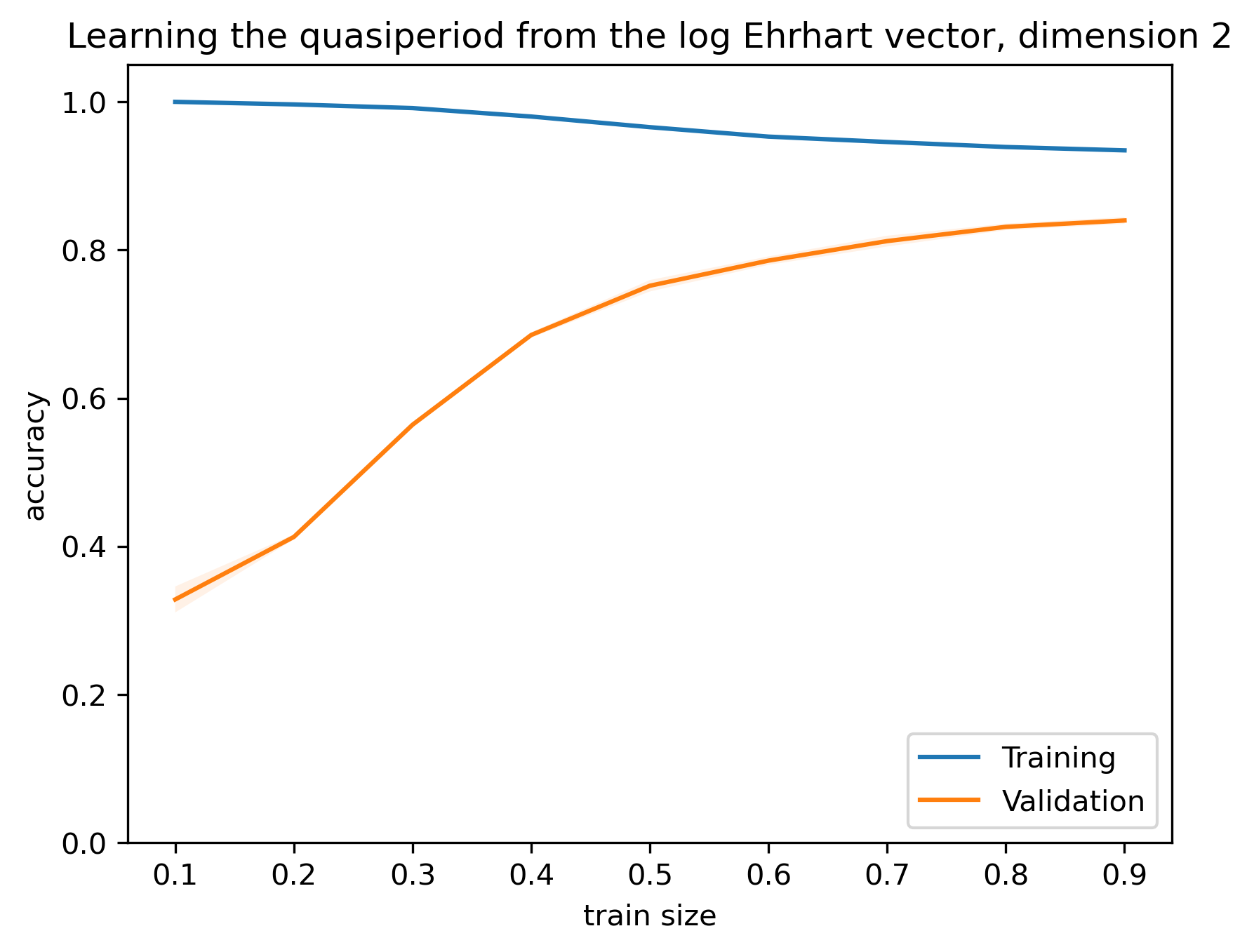}\\
	\includegraphics[height=5.5cm]{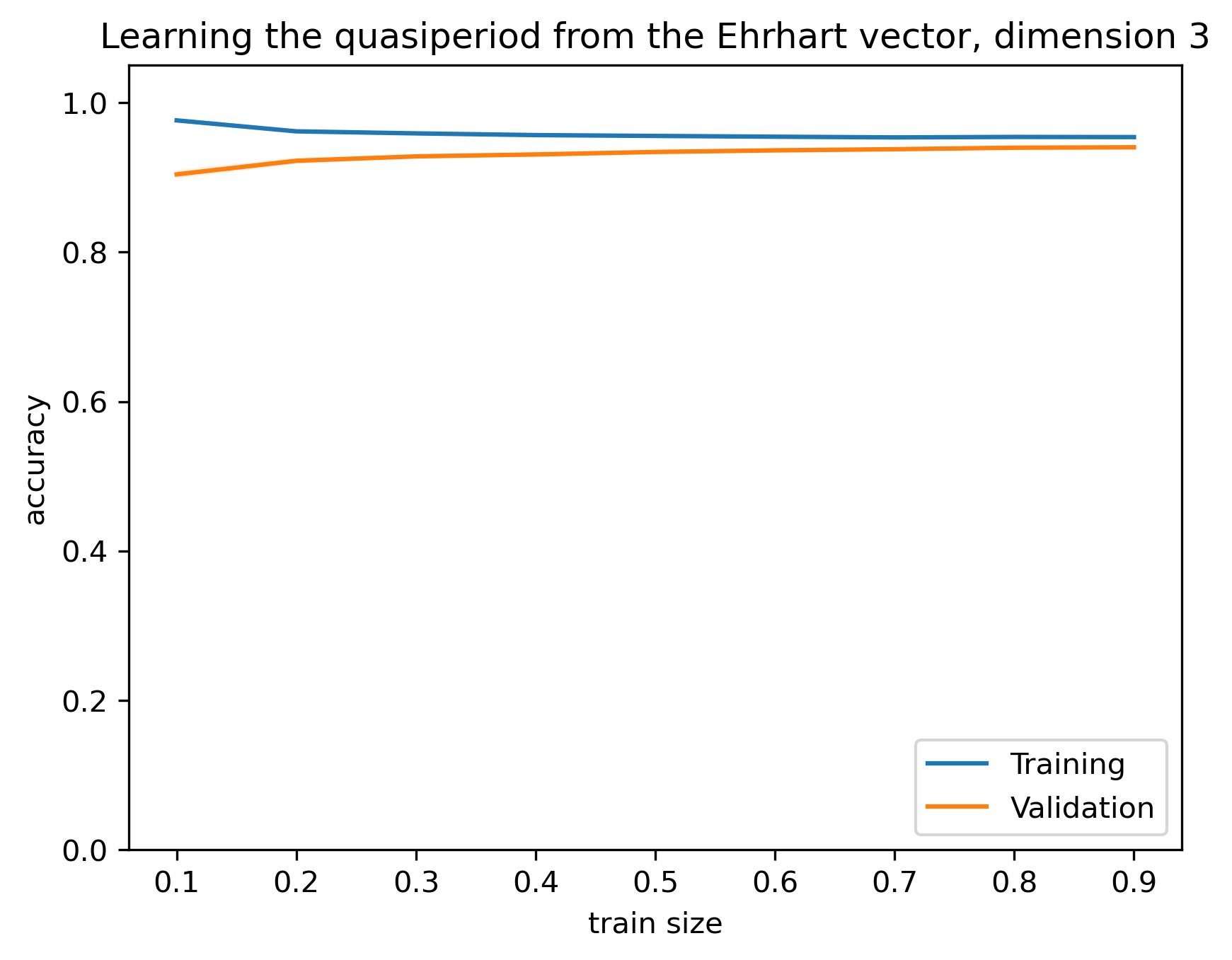}
	\includegraphics[height=5.5cm]{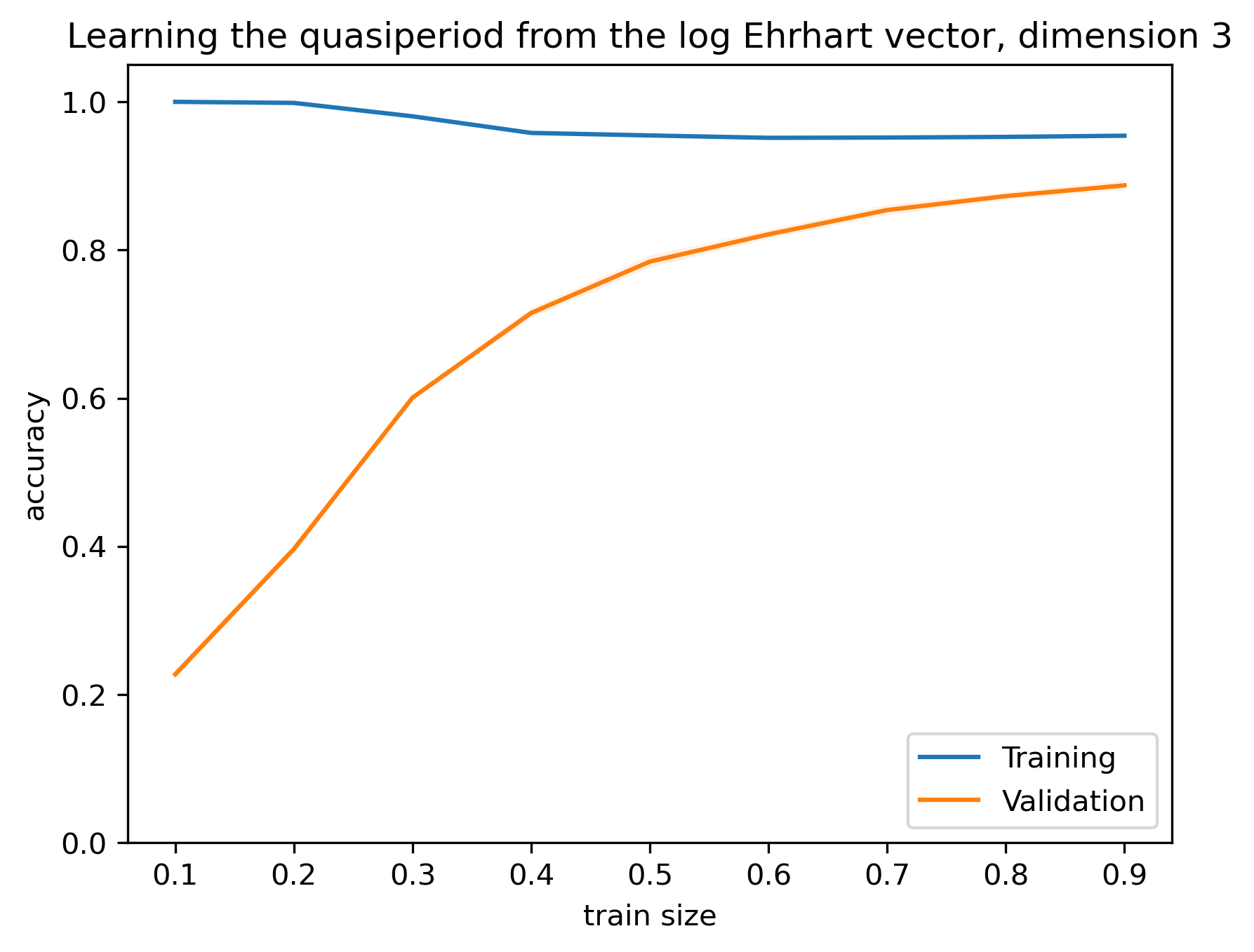}\\
	\includegraphics[height=5.5cm]{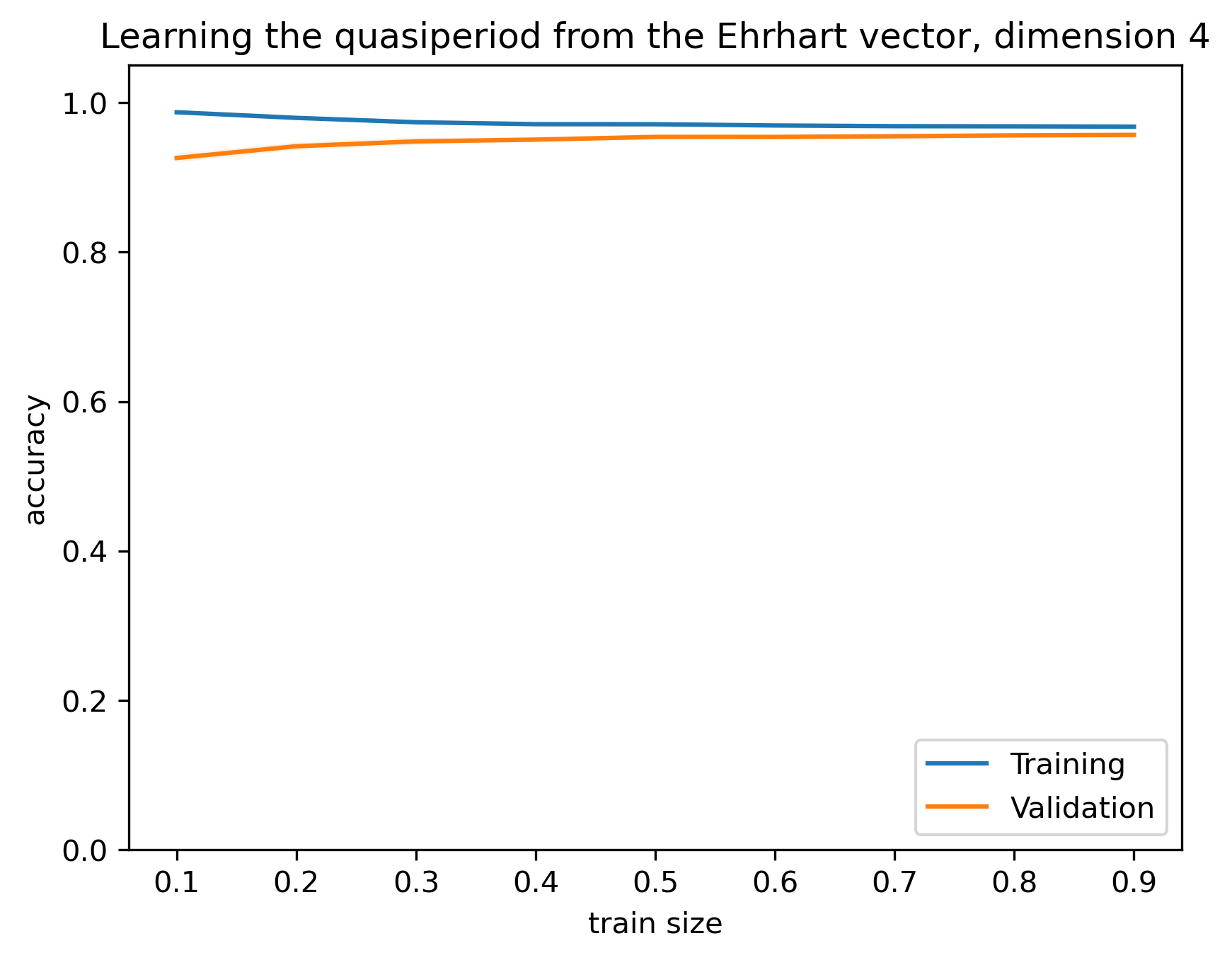}
	\includegraphics[height=5.5cm]{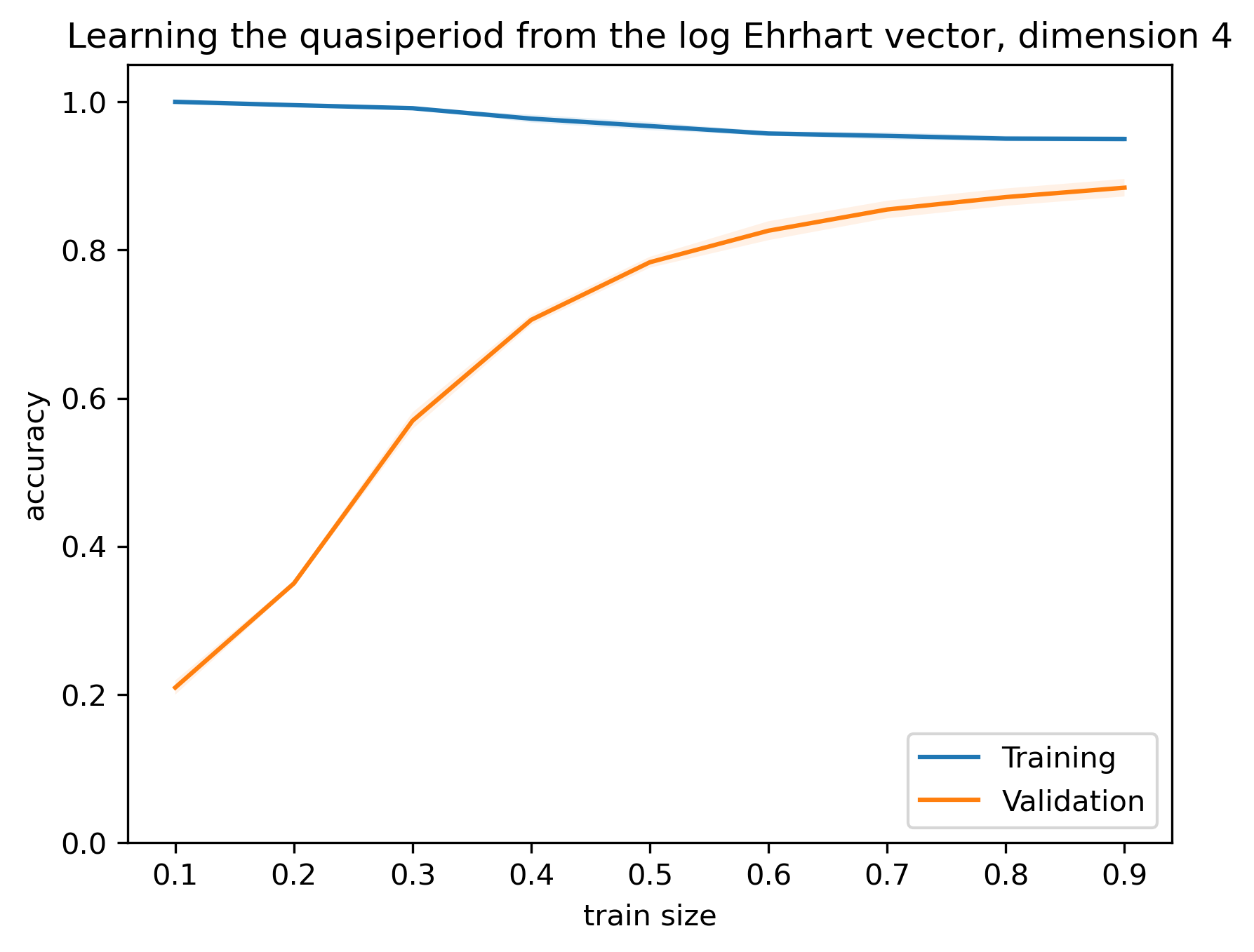}\\
	\caption{Quasi-period learning curves for a linear~SVM classifier, for both the Ehrhart vector (left-hand column) and the logarithmic Ehrhart vector (right-hand column).\label{fig:learning_curves}}
\end{figure}

\begin{figure}[tb]
	\includegraphics[width=.45\textwidth]{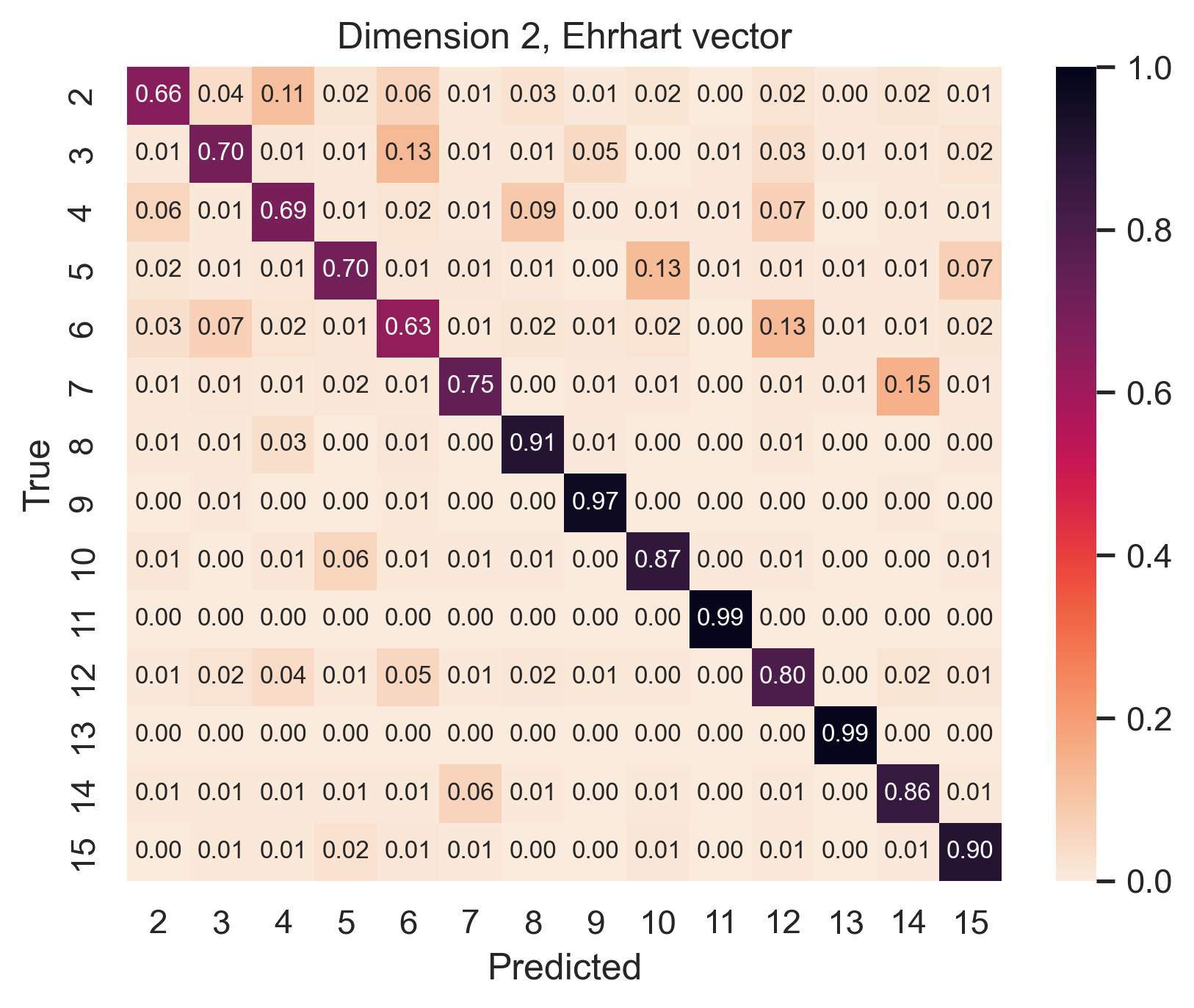}
	\includegraphics[width=.45\textwidth]{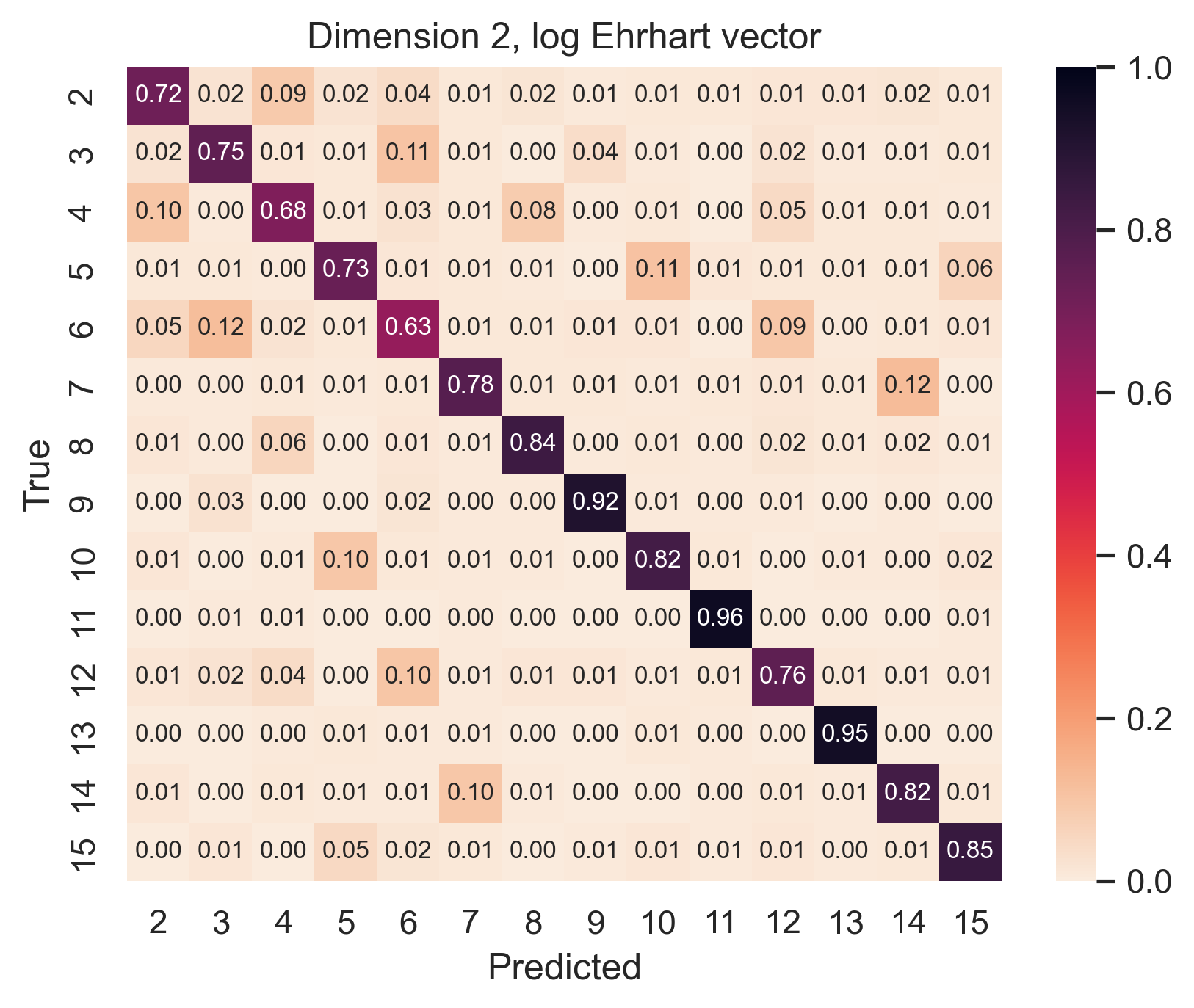}\\
	\includegraphics[width=.45\textwidth]{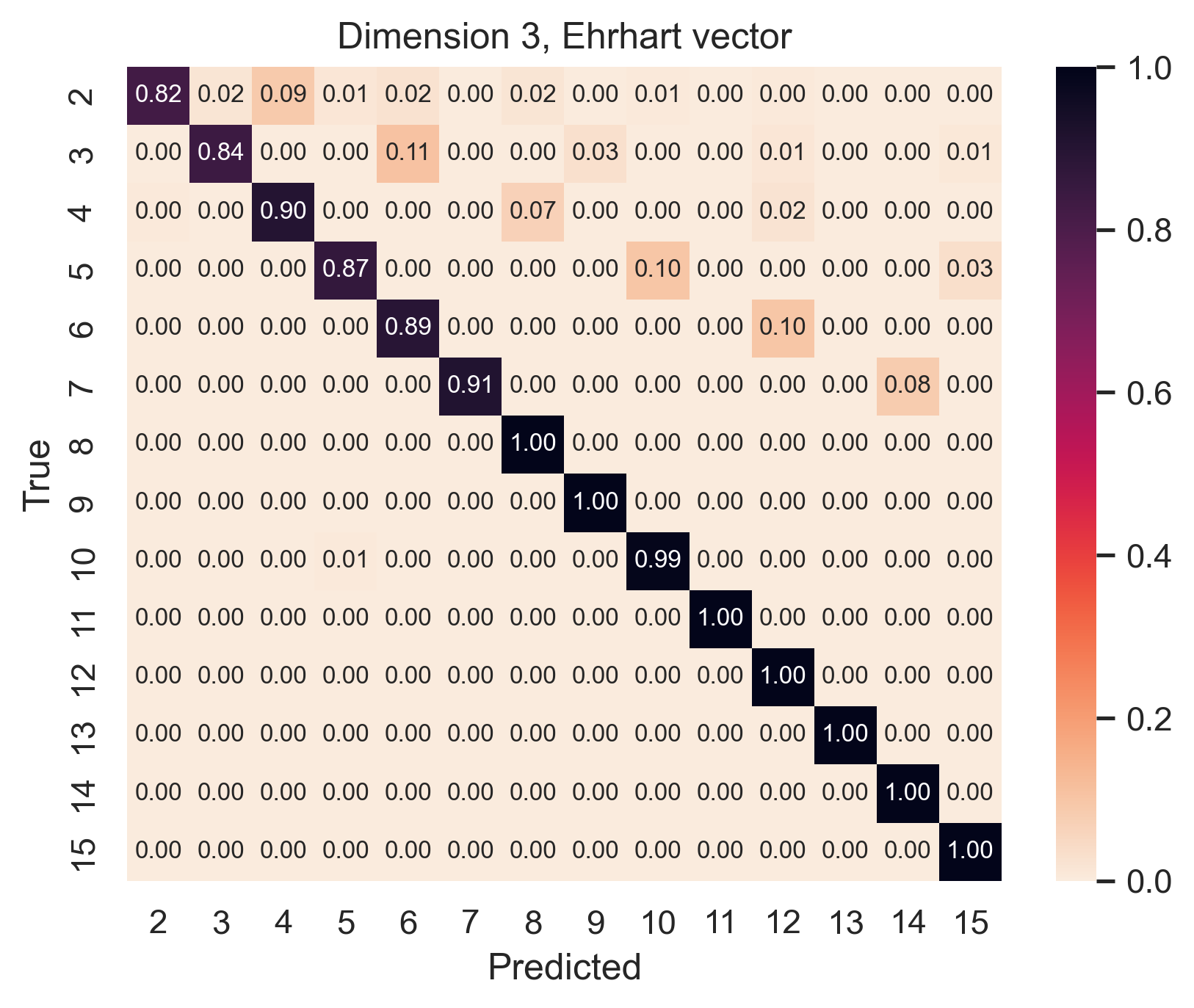}
	\includegraphics[width=.45\textwidth]{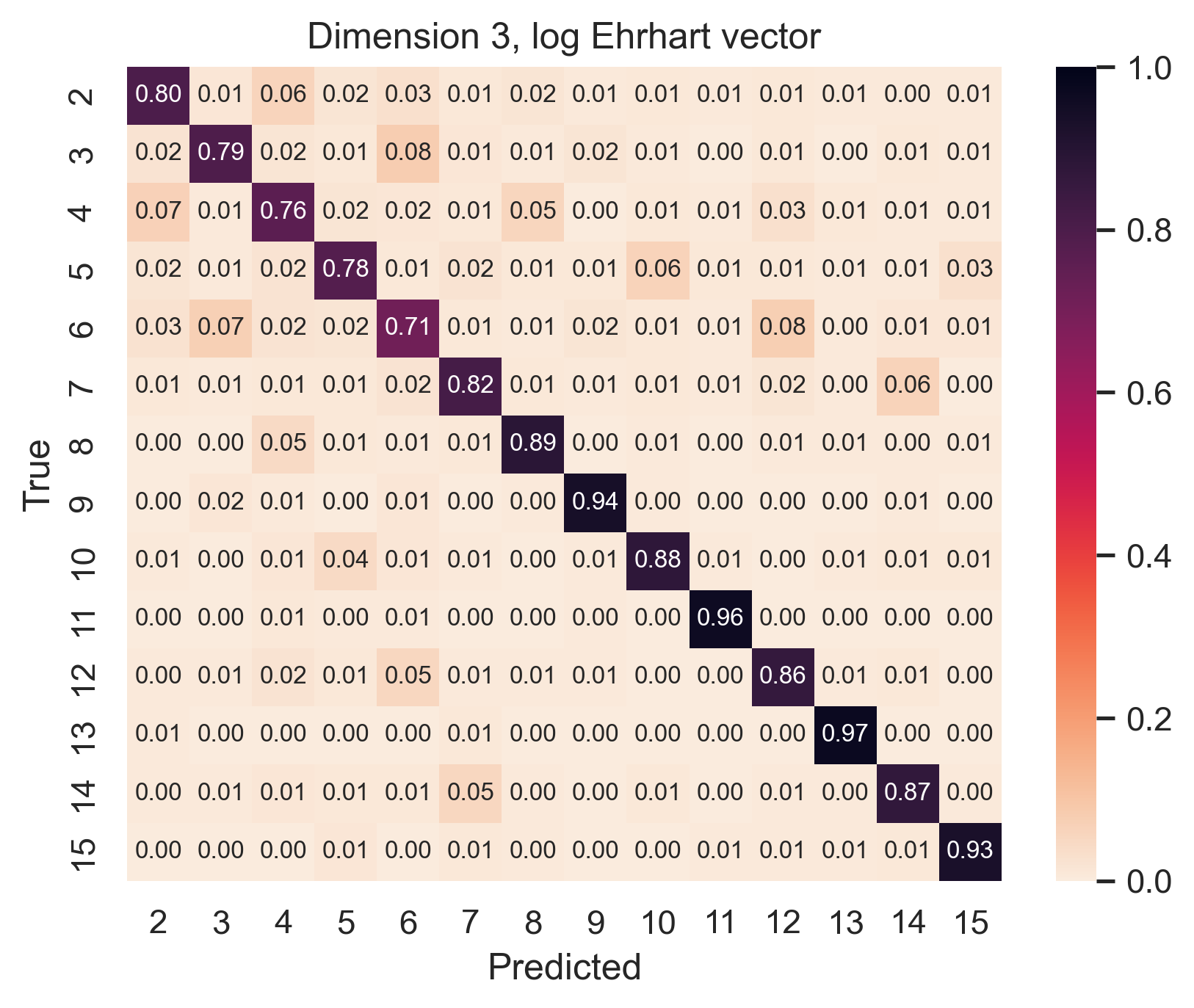}\\
	\includegraphics[width=.45\textwidth]{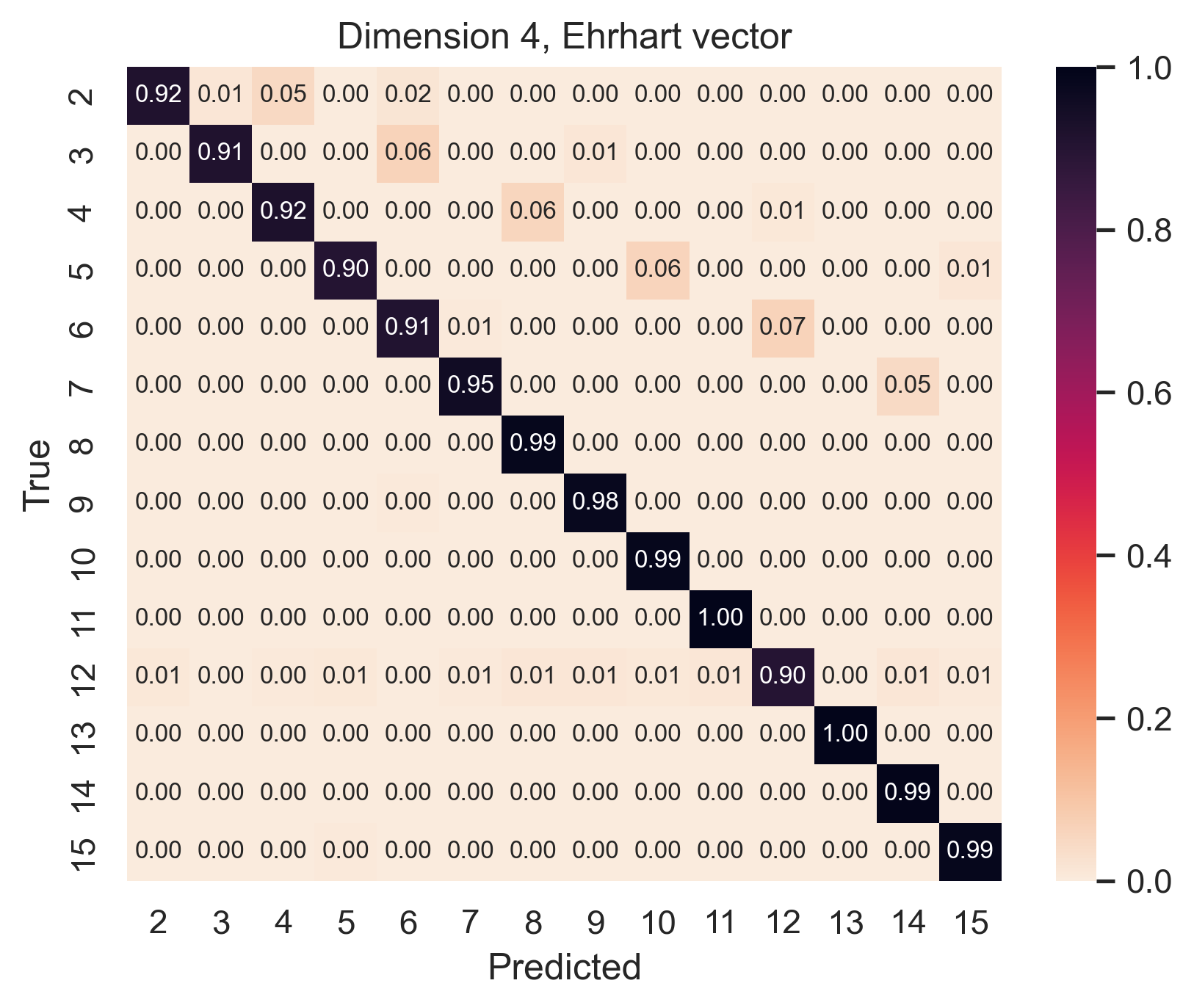}
	\includegraphics[width=.45\textwidth]{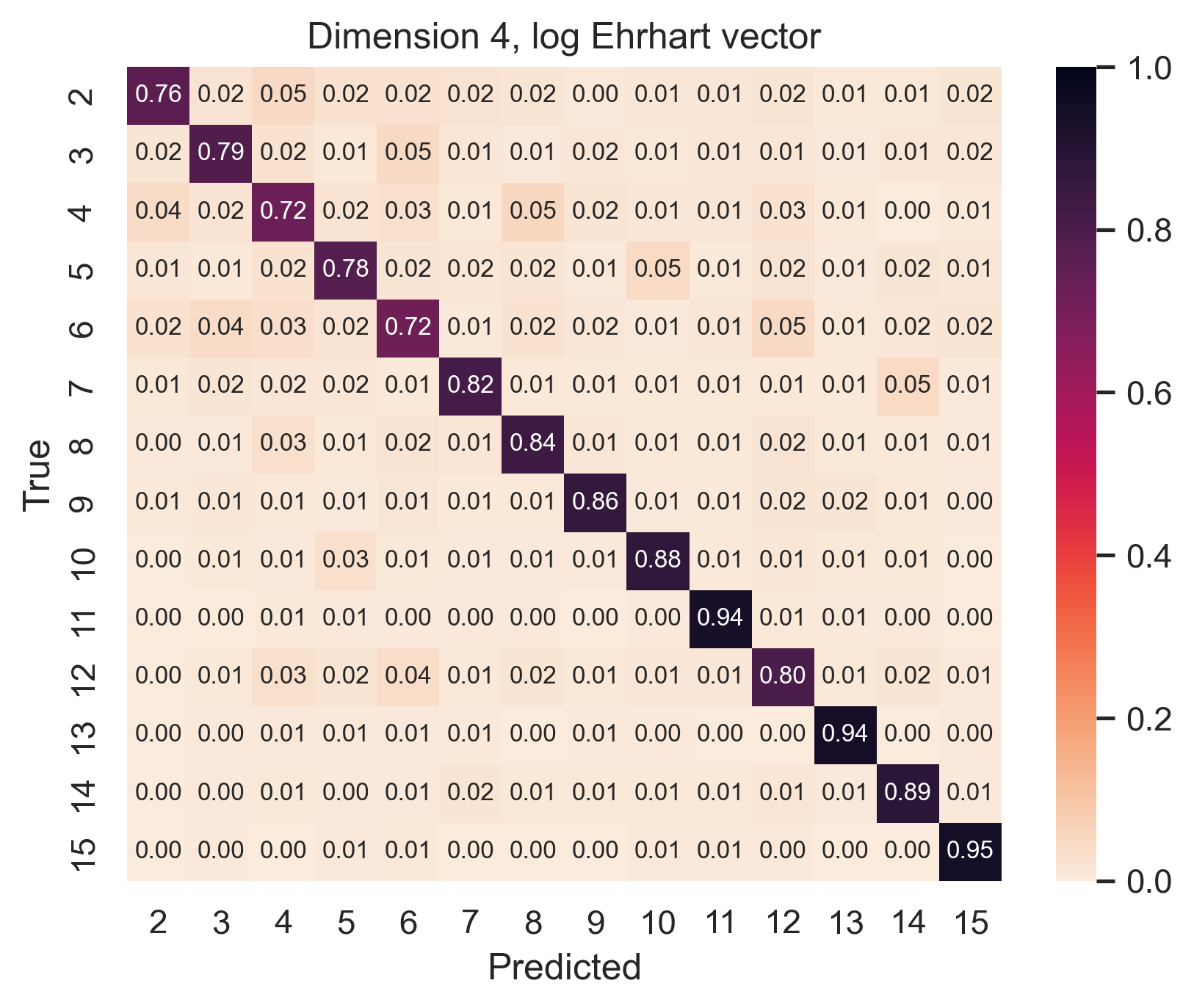}\\
	\caption{Row-normalised confusion matrices for a linear~SVM classifier for the quasi-period, for both the Ehrhart vector (left-hand column) and the logarithmic Ehrhart vector (right-hand column).\label{fig:confusion_matrices}}
\end{figure}

%-------------------------------------------------------------------------------
To learn the quasi-period of a rational polytope from its Ehrhart vector, we used a scikit-learn pipeline consisting of a \texttt{StandardScaler} followed by a \texttt{LinearSVC} classifier. We fixed a dimension~$d \in \{2,3,4\}$, moved to PCA co-ordinates, and used~$50\%$ of the data ($N=14000$) for training the classifier and hyperparameter tuning, holding out the remaining~$50\%$ of the data for model validation. Results are summarised on the left-hand side of Table~\ref{tab:quasiperiod}, with learning curves in the left-hand column of Figure~\ref{fig:learning_curves} and confusion matrices in the left-hand column of Figure~\ref{fig:confusion_matrices}. The confusion matrices hint at some structure in the misclassified data. 

One could also use the same pipeline but applied to the logarithmic Ehrhart vector rather than the Ehrhart vector. Results are summarised on the right-hand side of Table~\ref{tab:quasiperiod}, with learning curves in the right-hand column of Figure~\ref{fig:learning_curves} and confusion matrices in the right-hand column of Figure~\ref{fig:confusion_matrices}. Using the logarithmic Ehrhart data resulted in a less accurate classifier, but the learning curves suggest that this might be improved by adding more training data. Again there are hints of structure in the misclassified data.

\begin{table}[tb]
	\small
	\centering
	\begin{tabular}{ccccccc}
	\toprule
	\multicolumn{3}{c}{Ehrhart vector} & \phantom{xxx} &
	\multicolumn{3}{c}{log Ehrhart vector} \\ 
	\cmidrule{1-3} \cmidrule{5-7}
	Dimension& $C$ & Accuracy &&
	Dimension& $C$ & Accuracy \\
	\cmidrule{1-3} \cmidrule{5-7}
	2 & 0.01  & 80.6\% && 2 & 1 & 79.7\% \\
	3 & 0.001 & 94.0\% && 3 & 1 & 85.2\% \\
	4 & 0.001 & 95.3\% && 4 & 1 & 83.2\% \\
		\bottomrule
		\end{tabular}
	\caption{The regularization hyperparameter~$C$ and accuracy for a \texttt{LinearSVC} classifier predicting the quasi-period from the Ehrhart vector and logarithmic Ehrhart vector of rational polytopes~$P$.}
	\label{tab:quasiperiod}
\end{table}

%-------------------------------------------------------------------------------
\subsection{Forward differences}\label{sec:forward_differences}
%-------------------------------------------------------------------------------
Let~$y$ denote the sequence~$(y_m)_{m=0}^\infty$. Recall the forward difference operator~$\Delta$ defined on the space of sequences:
\[
	\Delta y = (y_{m+1} - y_m)_{m=0}^\infty.
\]
A sequence~$y$ depends polynomially on~$m$, that is
\begin{align*}
	y_m = a_0 + a_1 m + \cdots + a_d m^d &&
	\text{for some~$a_0,\ldots,a_d \in \RR$,}
\end{align*}
if and only if~$y$ lies in the kernel of~$\Delta^{d+1}$. Furthermore, in this case,~$\Delta^d y$ is the constant sequence with value~$d! a_d$.

Thus a sequence~$(y_m)_{m=0}^\infty$ is quasi-polynomial of degree~$d$ and period~$k$, in the sense of~\S\ref{sec:quasiperiod}, if and only if it lies in the kernel of~$\Delta_k^{d+1}$ where~$\Delta_k$ is the~$k$-step forward difference operator,
\[
	\Delta_k \left( (y_m)_{m=0}^\infty \right) = (y_{m+k} - y_m)_{m=0}^\infty,
\]
and~$(y_m)_{m=0}^\infty$ does not lie in the kernel of~$\Delta_k^d$. Furthermore, in this case, we can determine the leading coefficients of the polynomials~$f_0,\ldots,f_{k-1}$ by examining the values of the~$k$-periodic sequence~$\Delta_k^{d} y$. When~$y$ arises as the Ehrhart series of a rational polytope~$P$, all of these constant terms equal the volume of~$kP$, and so the value of the constant sequence~$\Delta_k^d y$ determines the normalised volume~$\Vol(P)$.

This discussion suggests that an SVM classifier with linear kernel should be able to learn the quasi-period and volume with high accuracy from the Ehrhart vector of a rational polytope, at least if we consider only polytopes of a fixed dimension~$d$. Having quasi-period~$k$ amounts to the Ehrhart vector lying in (a relatively open subset of) a certain subspace~$\ker \Delta_k^{d+1}$; these subspaces, being linear objects, should be easily separable using hyperplanes. Similarly, having fixed normalised volume amounts to lying in a given affine subspace; such affine subspaces should be easily separable using affine hyperplanes.

From this point of view it is interesting that an SVM classifier with linear kernel also learns the quasi-period with reasonably high accuracy from the logarithmic Ehrhart vector. Passing from the Ehrhart vector to the logarithmic Ehrhart vector replaces the linear subspaces~$\ker \Delta_k^{d+1}$ by non-linear submanifolds. But our experiments above suggest that these non-linear submanifolds must nonetheless be close to being separable by appropriate collections of affine hyperplanes.
%-------------------------------------------------------------------------------
\section{A remark on the Gorenstein index}\label{sec:gorenstein}
%-------------------------------------------------------------------------------
In this section we discuss a geometric question where machine learning techniques failed. This involves a more subtle combinatorial invariant called the Gorenstein index, and was the question that motivated the rest of the work in this paper. We then suggest why in retrospect we should not have expected to be able to answer this question using machine learning (or at all).
%-------------------------------------------------------------------------------
\subsection{The Gorenstein index}
%-------------------------------------------------------------------------------
Fix a rank~$d$ lattice~$N\cong\ZZ^d$, write~$\NQ \coloneq N\otimes_\ZZ\QQ$, and let~$P\subset \NQ$ be a lattice polytope. The~\emph{polar polyhedron} of~$P$ is given by
\[
P^*\coloneq\{u\in M\otimes_\ZZ\QQ\mid u(v)\geq -1\text{ for all }v\in P\}
\]
where~$M\coloneq\Hom(N,\ZZ)\cong\ZZ^d$ is the lattice dual to~$N$. The polar polyhedron~$P^*$ is a convex polytope with rational vertices if and only if the origin lies in the strict interior of~$P$. In this case,~$(P^*)^*=P$. The smallest positive integer~$k_P$ such that~$k_PP^*$ is a lattice polytope is called the~\emph{Gorenstein index} of~$P$.

The Gorenstein index arises naturally in the context of~\emph{Fano toric varieties}, and we will restrict our discussion to this setting. Let~$P\subset\NQ$ be a lattice polytope such that the vertices of~$P$ are primitive lattice vectors and that the origin lies in the strict interior of~$P$; such polytopes are called~\emph{Fano}~\cite{KN13}. The spanning fan~$\Sigma_P$ of~$P$ -- that is, the complete fan whose cones are generated by the faces of~$P$ -- gives rise to a Fano toric variety~$X_P$~\cite{Ful93,CLS11}. This construction gives a one-to-one correspondence between~$\GL_d(\ZZ)$-equivalence classes of Fano polytopes and isomorphism classes of Fano toric varieties. Let~$P\subset\NQ$ be a Fano polytope that corresponds to a Fano toric variety~$X\coloneq X_P$. The polar polytope~$P^*\subset\MQ$ then corresponds to a divisor on~$X$ called the~\emph{anticanonical divisor}, which is denoted by~$-K_X$. In general~$-K_X$ is an ample~$\QQ$-Cartier divisor, and~$X$ is~$\QQ$-Gorenstein. The Gorenstein index~$k_P$ of~$P$ is equal to the smallest positive multiple~$k$ of the anticanonical divisor such that~$-kK_X$ is Cartier.

Under the correspondence just discussed, the Ehrhart series~$\Ehr_{P^*}$ of the polar polytope~$P^*$ coincides with the Hilbert series~$\Hilb_X(-K_X)$. The Hilbert series is an important numerical invariant of~$X$, and it makes sense to ask whether the Gorenstein index of~$X$ is determined by the Hilbert series. Put differently:

\begin{question}\label{qu:ML_Gorenstein}
    Given a Fano polytope~$P$, can~ML recover the Gorenstein index~$k_P$ of~$P$ from sufficiently many terms of the Ehrhart series~$\Ehr_{P^*}$ of the polar polytope~$P^*$?
\end{question}

There are good reasons, as we discuss below, to expect the answer to Question~\ref{qu:ML_Gorenstein} to be `no'. But part of the power of~ML in mathematics is that it can detect or suggest structure that was not known or expected previously (see e.g.~\cite{DaviesEtAl2021}). That did not happen on this occasion: applying the techniques discussed in \S\ref{sec:ML_dimension} and \S\ref{sec:ML_quasiperiod} did not allow us to predict the Gorenstein index of~$P$ from the Ehrhart series of the polar polytope~$\Ehr_{P^*}$.
%-------------------------------------------------------------------------------
\subsection{Should we have expected this?}\label{sec:Gorenstein_discussion}
%-------------------------------------------------------------------------------
The Hilbert series is preserved under an important class of deformations called~\emph{qG-deformations}~\cite{KS-B88}. But the process of~\emph{mutation}~\cite{ACGK12} can transform a Fano polytope~$P$ to a Fano polytope~$Q$ with~$\Ehr_{P^*}=\Ehr_{Q^*}$. Mutation gives rise to a qG-deformation from~$X_P$ to~$X_Q$~\cite{ACC+16}, but need not preserve the Gorenstein index:~$k_P$ need not be equal to~$k_Q$. Thus the Gorenstein index is not invariant under qG-deformation. It might have been unrealistic to expect that a qG-deformation invariant quantity (the Hilbert series) could determine an invariant (the Gorenstein index) which can vary under qG-deformation.
%-------------------------------------------------------------------------------
\subsection*{Quasiperiod collapse}
%-------------------------------------------------------------------------------
Although the phenomenon of quasi-period collapse remains largely mysterious from a combinatorial view-point, in the context of toric geometry one possible explanation arises from mutation and qG-deformation~\cite{KW18}. The following example revisits Examples~\ref{eg:P2} and~\ref{eg:P114} from this point of view, and illustrates why Question~\ref{qu:ML_Gorenstein} cannot have a meaningful positive answer.

\begin{figure}[tb]
    \includegraphics[scale=0.8]{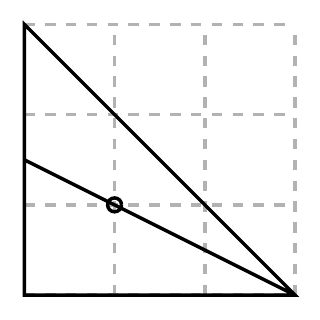}
    \raisebox{23.5px}{$\longrightarrow$}
    \includegraphics[scale=0.8]{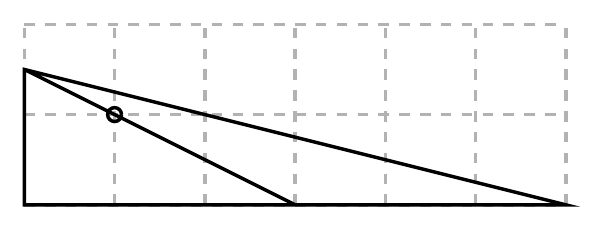}
	\caption{The mutation from~$P_{(1,1,1)}$ to~$P_{(1,1,2)}$ gives a scissors congruence between the polar polytopes.}
    \label{fig:scissors_congruence}
\end{figure}

\begin{example}\label{eg:mutation_P2}
Let~$P_{(a,b,c)}\subset\NQ$ denote the~$2$-dimensional Fano polytope associated with weighted projective space~$\PP(a^2,b^2,c^2)$, where~$a$,~$b$,~$c$ are pairwise coprime positive integers. Then~$P_{(1,1,1)}\coloneq\conv\{(1,0),(0,1),(-1,-1)\}$ is the Fano polygon associated with~$\PP^2$, with polar polygon~$P^*_{(1,1,1)}\subset\MQ$ the lattice triangle appearing in Example~\ref{eg:P2}.

The graph of mutations of~$P_{(1,1,1)}$ has been completely described~\cite{HP10,AK16}. Up to~$\GL_2(\ZZ)$-equivalence, there is exactly one mutation from~$P_{(1,1,1)}$: this gives~$P_{(1,1,2)}$ corresponding to~$\PP(1,1,4)$. The polar polygon~$P^*_{(1,1,2)}$ is the rational triangle in Example~\ref{eg:P114}. As Figure~\ref{fig:scissors_congruence} illustrates, a mutation between polytopes gives rise to a scissors congruence~\cite{HaaseMcAllister2008} between polar polytopes. Mutation therefore preserves the Ehrhart series of the polar polytope, and this explains why we have quasi-period collapse in this example. 

We can mutate~$P_{(1,1,2)}$ in two ways that are distinct up to the action of~$\GL_2(\ZZ)$: one returns us to~$P_{(1,1,1)}$, whilst the other gives~$P_{(1,2,5)}$. Continuing to mutate, we obtain an infinite graph of triangles~$P_{(a,b,c)}$, where the~$(a,b,c)$ are the~\emph{Markov triples}, that is, the positive integral solutions to the Markov equation:
\[
3xyz=x^2+y^2+z^2.
\]
The Gorenstein index of~$P_{(a,b,c)}$ is~$abc$. In particular, the Gorenstein index can be made arbitrarily large whilst the Ehrhart series, and hence quasi-period, of~$P^*_{(a,b,c)}$ is fixed. See~\cite{KW18} for details.
\end{example}
%-------------------------------------------------------------------------------
\section{Conclusion}\label{sec:conclusion}
%-------------------------------------------------------------------------------
We have seen that Support Vector Machine methods are very effective at extracting the dimension and volume of a lattice or rational polytope~$P$, and the quasi-period of a rational polytope~$P$, from the initial terms of its Ehrhart series. We have also seen that~ML methods are unable to reliably determine the Gorenstein index of a Fano polytope~$P$ from the Ehrhart series of its polar polytope~$P^*$. The discussions in~\S\ref{sec:log_asymptotics},~\S\ref{sec:forward_differences}, and~\S\ref{sec:Gorenstein_discussion} suggest that these results are as expected: that~ML is detecting known and understood structure in the dimension, volume, and quasiperiod cases, and that there is probably no structure to detect in the Gorenstein index case. But there is a more useful higher-level conclusion to draw here too: when applying~ML methods to questions in pure mathematics, one needs to think carefully about methods and results. Questions~\ref{qu:ML_quasiperiod} and~\ref{qu:ML_Gorenstein} are superficially similar, yet one is amenable to~ML and the other is not. Furthermore, applying standard~ML recipes in a naive way would have led to false negative results. For example, since the Ehrhart series grows so fast, it would have been typical to suppress the growth rate by taking logarithms, and also to pass to principal components. Taking logarithms is a good idea for some of our questions but not for others; this reflects the mathematical realities underneath the data, and not just whether the vector-components~$y_m$ involved grow rapidly with~$m$ or not. Passing to principal components is certainly a useful tool, but naive feature extraction would have retained only the first principal component, which is responsible for more than~$99.999\%$ of the variation (in both the logarithmic Ehrhart vector and the Ehrhart vector). This would have left us unable to detect the positive answers to Questions~\ref{qu:ML_dimension} and~\ref{qu:ML_quasiperiod}, in the former case because projection to one dimension amalgamates clusters (see Figure~\ref{fig:pca_dim}) and in the latter case because we need to detect whether the Ehrhart vector lies in a certain high-codimension linear subspace and that structure is destroyed by projection to a low-dimensional space.
%-------------------------------------------------------------------------------
\subsection*{Acknowledgments}\label{sec:acknowledgements}
%-------------------------------------------------------------------------------
TC~is supported by ERC Consolidator Grant~682603 and EPSRC Programme Grant~EP/N03189X/1. JH~is supported by a Nottingham Research Fellowship. AK~is supported by EPSRC Fellowship~EP/N022513/1. 
%-------------------------------------------------------------------------------
\providecommand{\bysame}{\leavevmode\hbox to3em{\hrulefill}\thinspace}
\providecommand{\MR}{\relax\ifhmode\unskip\space\fi MR }
\providecommand{\MRhref}[2]{%
  \href{http://www.ams.org/mathscinet-getitem?mr=#1}{#2}
}
\providecommand{\href}[2]{#2}

%-------------------------------------------------------------------------------
\end{document}